\newif\ifpdf
\ifx\pdfoutput\undefined
        \pdffalse
\else
        \pdfoutput=1 
        \pdftrue
\fi

\ifpdf
\documentclass[pdftex]{amsart}
\else
\documentclass{amsart}
\fi

\usepackage{amssymb}

\ifpdf
        \usepackage{color}
        \usepackage[pdftex]{graphicx}
\else
        \usepackage{graphicx}
\fi

\usepackage[all]{xy}

\title[Parall{\'e}lisme et homotopie (II)]{Concurrent Process up to Homotopy (II)}
\author[P. Gaucher]{Philippe Gaucher}
\address{Institut de Recherche Math\'ematique Avanc\'ee\\ ULP et
CNRS\\ 7 rue Ren\'e Descar\-tes\\ 67084 Strasbourg Cedex\\ France}
\email{gaucher@math.u-strasbg.fr}
\urladdr{http://www-irma.u-strasbg.fr/\~{}gaucher/}

\newcommand{\p}\times
\renewcommand{\vec}{\overrightarrow}
\renewcommand{\P}{\mathbb{P}}

\newcommand{\be}{\begin{equation}}
\newcommand{\ee}{\end{equation}}
\newcommand{\bea}{\begin{eqnarray}}
\newcommand{\eea}{\end{eqnarray}}
\newcommand{\beas}{\begin{eqnarray*}}
\newcommand{\eeas}{\end{eqnarray*}}

\newtheorem{theoreme}{Th{\'e}or{\`e}me}[section]
\newtheorem{prop}[theoreme]{Proposition}

\newtheorem{question}[theoreme]{Question}

\newtheorem{rem}[theoreme]{Remarque}

\newtheorem{definition}[theoreme]{D{\'e}finition}

\newcommand{\bd}{\begin{defn}}
\newcommand{\ed}{\end{defn}}
\newcommand{\bcd}{\begin{defn}}
\newcommand{\ecd}{\end{defn}}
\newcommand{\bex}{\begin{exmp}}
\newcommand{\eex}{\end{exmp}}
\newcommand{\bp}{\begin{prop}}
\newcommand{\ep}{\end{prop}}
\newcommand{\bth}{\begin{theoreme}}
\renewcommand{\eth}{\end{theoreme}}
\newcommand{\br}{\begin{rem}}
\newcommand{\er}{\end{rem}}
\newcommand{\bpf}{\begin{proof}}
\newcommand{\epf}{\end{proof}}

\newcommand{\fl}[1]{\ar@{->}[l]_{#1}}
\newcommand{\fr}[1]{\ar@{->}[r]^-{#1}}
\newcommand{\fd}[1]{\ar@{->}[d]_{#1}}
\newcommand{\fu}[1]{\ar@{->}[u]^{#1}}
\newcommand{\f}[2]{\ar@{->}[#1]|{#2}}
\newcommand{\ff}[2]{\ar@2{->}[#1]|{#2}}
\newcommand{\frr}[1]{\ar@{->}[rr]^{#1}}

\newcommand{\iso}{\cong}

\newcommand{\vI}{\vec{I}}

\renewcommand{\leq}{\leqslant}
\renewcommand{\geq}{\geqslant}

\def\cartesien{%
  \ar@{-}[]+R+<6pt,-2pt>;[]+RD+<6pt,-6pt>%
  \ar@{-}[]+D+<2pt,-6pt>;[]+RD+<6pt,-6pt>%
}
\def\cocartesien{%
  \ar@{-}[]+L+<-6pt,+2pt>;[]+LU+<-6pt,+6pt>%
  \ar@{-}[]+U+<-2pt,+6pt>;[]+LU+<-6pt,+6pt>%
}

\newcommand{\brm}[1]{\rm{\mathbf{#1}}}

\newcommand{\dtop}{{\brm{Flow}}}

\newcommand{\tdtop}{{\brm{FLOW}}}
\newcommand{\ttop}{{\brm{TOP}}}

\makeatletter

\def\varholim@#1#2{%
  \vtop{\m@th\ialign{##\cr
    \hfil$#1\operator@font holim$\hfil\cr
    \noalign{\nointerlineskip\kern1.5\ex@}#2\cr
    \noalign{\nointerlineskip\kern-\ex@}\cr}}%
}
\def\holimproj{%
  \mathop{\mathpalette\varholim@{\leftarrowfill@\textstyle}}\nmlimits@
}
\def\holimind{%
  \mathop{\mathpalette\varholim@{\rightarrowfill@\textstyle}}\nmlimits@
}

\newskip\@bigflushglue \@bigflushglue = -100pt plus 1fil

\def\bigcentering{\let\\\@centercr\rightskip\@bigflushglue%
\leftskip\@bigflushglue
\parindent\z@\parfillskip\z@skip}

\makeatother

\DeclareMathOperator{\ctop}{{\textbf{Top}}}

\begin{document}

\begin{abstract}
On d{\'e}montre que la cat{\'e}gorie des CW-complexes globulaires {\`a}
dihomotopie pr{\`e}s est {\'e}quivalente {\`a} la cat{\'e}gorie des flots {\`a}
dihomotopie faible pr{\`e}s. Ce th{\'e}or{\`e}me est une g{\'e}n{\'e}ralisation du
th{\'e}or{\`e}me classique disant que la cat{\'e}gorie des CW-complexes modulo
homotopie est {\'e}quivalente {\`a} la cat{\'e}gorie des espaces topologiques
modulo homotopie faible.

One proves that the category of globular CW-complexes up to dihomotopy
is equivalent to the category of flows up to weak dihomotopy. This
theorem generalizes the classical theorem which states that the category
of CW-complexes up to homotopy is equivalent to the category of
topological spaces up to weak homotopy.
\end{abstract}
\maketitle

\section{Rappels sur les flots}

Cette note est la deuxi{\`e}me de deux notes pr{\'e}sentant quelques r{\'e}sultats de
\cite{flow}.  Tous les espaces topologiques sont suppos{\'e}s faiblement 
s{\'e}par{\'e}s et compactement engendr{\'e}s, c'est-{\`a}-dire dans ce cas
hom{\'e}omorphes {\`a} la limite inductive de leurs sous-espaces compacts
(cf. l'appendice de \cite{Ref_wH} pour un survol des propri{\'e}t{\'e}s de ces
espaces). La cat{\'e}gorie correspondante est not{\'e}e $\ctop$.  On travaille
ainsi dans une cat{\'e}gorie d'espaces topologiques qui est non seulement
compl{\`e}te et cocompl{\`e}te mais en plus cart{\'e}siennement ferm{\'e}e
\cite{MR35:970}. En d'autres termes, le foncteur $-\p
X:\ctop\rightarrow \ctop$ a un adjoint {\`a} droite $\ttop(X,-)$.  Dans la
suite, $D^n$ est le disque ferm{\'e} de dimension $n$ et $S^{n-1}$ est le
bord de $D^n$, {\`a} savoir la sph{\`e}re de dimension $n-1$. En particulier
la sph{\`e}re de dimension $0$ est la paire $\{-1,+1\}$.

\begin{definition}\cite{flow} 
Un flot $X$ consiste en la donn{\'e}e d'un ensemble $X^0$ appel{\'e}
$0$-squelette, d'un espace topologique $\P X$ appel{\'e} espace des
chemins, de deux applications continues $s:\P X\rightarrow X^0$ et
$t:\P X\rightarrow X^0$ ($X^0$ {\'e}tant muni de la topologie discr{\`e}te) et
d'une application continue $*:\{(x,y)\in \P X\p\P
X,t(x)=s(y)\}\rightarrow \P X$ satisfaisant les axiomes $s(x*y)=s(x)$,
$t(x*y)=t(y)$ et enfin $x*(y*z)=(x*y)*z$ pour tout $x,y,z\in \P X$. Un
morphisme de flots $f$ de $X$ vers $Y$ est une application continue de
$X^0\sqcup \P X$ vers $Y^0\sqcup \P Y$ telle que $f(X^0)\subset Y^0$,
$f(\P X)\subset \P Y$, $s(f(x))=f(s(x))$, $t(f(x))=f(t(x))$,
$f(x*y)=f(x)*f(y)$. La cat{\'e}gorie correspondante est not{\'e}e $\dtop$.
\end{definition}

Les {\'e}l{\'e}ments de $\P X$ sont appel{\'e}s les chemins d'ex{\'e}cution (non-constants) et 
ceux de $X^0$ les {\'e}tats de $X$ (ou encore les chemins d'ex{\'e}cution constants).

Si $Z$ est un espace topologique, le flot $Glob(Z)$ est d{\'e}fini comme
suit : $\P Glob(Z)=Z$, $Glob(Z)^0=\{0,1\}$, et enfin $s=0$ et $t=1$
(il n'y a pas de chemins d'ex{\'e}cution composables).  Si $Z$ est un
singleton, on obtient le flot $\vI$ correspondant au segment dirig{\'e}.

\begin{definition}\cite{flow}
Deux morphismes de flots $f$ et $g$ de $X$ dans $Y$ sont S-homotopes s'il existe
une application continue $H:X\p [0,1]\rightarrow Y$ telle que
\begin{enumerate}
\item $H(-,u)$ est un morphisme de flots de $X$ dans $Y$
\item $H(-,0)=f$ et $H(-,1)=g$
\end{enumerate}
On {\'e}crit $f\sim_S g$.
\end{definition}

De l{\`a} on d{\'e}finit la notion de flots S-homotopes :

\begin{definition}\cite{flow}
Deux flots $X$ et $Y$ sont S-homotopes s'il existe un morphisme de
flots
$f:X\rightarrow Y$ et un morphisme de flots  $g:Y\rightarrow X$ tels que
$f\circ g\sim_{S} Id_Y$ et $g\circ f\sim_{S} Id_X$. On dit alors que $f$ 
et $g$ sont deux {\'e}quivalences de S-homotopie r{\'e}ciproques. \end{definition}

Il existe une deuxi{\`e}me classe de morphismes de flots appel{\'e}s T-homotopie qui 
permet d'identifier des flots ayant les m{\^e}mes propri{\'e}t{\'e}s informatiques 
sous-jacentes \cite{note1,flow}. Nous n'en dirons rien de nouveau ici par 
rapport {\`a} la pr{\'e}c{\'e}dente note. Il est donc inutile de rappeler la 
d{\'e}finition pr{\'e}cise de cette classe de morphismes.

\section{Flot S-cofibrant}

Un flot \textit{S-cofibrant} est obtenu par d{\'e}finition comme suit. On
part d'un espace discret $X^0$ et on consid{\`e}re le flot correspondant
ayant $X^0$ comme $0$-squelette et l'ensemble vide $\emptyset$ comme
espace de chemins. On lui attache ensuite des $\vI$ pour obtenir un
premier flot $X^1_0$. On passe ensuite de $X^1_n$ {\`a} $X^1_{n+1}$ pour
$n\geq 0$ en attachant des $Glob(D^{n+1})$ le long de morphismes
d'attachement $Glob(S^{n})\rightarrow X^1_n$. En d'autres termes, on
choisit une famille de morphismes de flots $f_i:Glob(S^{n})\rightarrow
X^1_n$ pour $i\in I$ et $X^1_{n+1}$ est alors obtenu par le diagramme
cocart{\'e}sien dans la cat{\'e}gorie des flots (qui se trouve {\^e}tre
cocompl{\`e}te)
\[
\xymatrix{
\bigsqcup_{i\in I}Glob(S^{n})\fr{\bigsqcup_{i\in I}f_i}\fd{}& X^1_n \fd{}\\
\bigsqcup_{i\in I}Glob(D^{n+1})\fr{} & \cocartesien X^1_{n+1}}
\]

Si $Z$ est un CW-complexe, alors $Glob(Z)$ est un exemple de flot
S-cofibrant. Vu la similitude entre la construction des CW-complexes
globulaires introduits dans \cite{diCW,note1} et celle des flots
S-cofibrants, on n'aura pas de mal {\`a} se convaincre du :

\begin{theoreme}\cite{flow} Le plongement des CW-complexes globulaires dans les flots
induit une {\'e}quivalence de cat{\'e}gories entre les CW-complexes
globulaires {\`a} S-homotopie et T-homotopie pr{\`e}s et les flots
S-cofibrants {\`a} S-homotopie et T-homotopie pr{\`e}s. \end{theoreme}

La preuve de ce th{\'e}or{\`e}me repose essentiellement sur le fait que deux
CW-complexes globulaires sont S-homotopes (resp. T-homotopes) si et
seulement si les flots correspondants sont S-homotopes
(resp. T-homotopes) (cf. \cite{note1}).

\section{Dihomotopie faible}

\begin{definition}\cite{flow} 
Un morphisme de flots $f:X\rightarrow Y$ est une S-homotopie faible si
$f$ induit une bijection entre $X^0$ et $Y^0$ et une {\'e}quivalence
d'homotopie faible de $\P X$ vers $\P Y$.  On dit alors que $X$ et $Y$
sont faiblement S-homotopes. \end{definition}

\begin{theoreme}\cite{flow} Un morphisme de flots $f:X\rightarrow Y$ entre 
deux flots S-cofibrants $X$ et $Y$ est une {\'e}quivalence de S-homotopie faible 
si et seulement si c'est une {\'e}quivalence de S-homotopie.  \end{theoreme}

Si les deux flots S-cofibrants sont de la forme $Glob(Z)$ et $Glob(T)$
o{\`u} $Z$ et $T$ sont deux CW-complexes, cet {\'e}nonc{\'e} correspond au
th{\'e}or{\`e}me classique de Whitehead.

Pour le d{\'e}montrer, on commence par introduire pour tout espace
topologique $U$ et tout flot $X$ un flot $U\boxtimes X$ de telle fa{\c c}on
que l'on ait la bijection naturelle $\dtop([0,1]\boxtimes X,Y)\iso
\ctop([0,1],\tdtop(X,Y))$ o{\`u} $\tdtop(X,Y)$ est l'ensemble des
morphismes de flots de $X$ {\`a} $Y$ muni de la Kelleyfication de la
topologie induite par celle de $\ttop(X,Y)$. Mais le foncteur
$\tdtop(X,-):\dtop\rightarrow \ctop$ ne commute pas avec les limites
projectives. En effet, si $X$ est un singleton, alors
$\tdtop(X,Y)=Y^0$ est toujours discret mais la limite projective
d'espaces topologiques discrets peut {\^e}tre totalement discontinue sans
{\^e}tre discr{\`e}te. On doit donc proc{\'e}der autrement pour d{\'e}finir
correctement $U\boxtimes X$. Mais d{\`e}s que $U$ est connexe, la
bijection naturelle ci-dessus est vraie et $U\boxtimes X$ est alors le
flot libre engendr{\'e} par $X^0\sqcup (U\p X)$ quotient{\'e} par les
relations $s(u,x)=s(x)$, $t(u,x)=t(x)$ et $(u,x)*(u,y)=(u,x*y)$ d{\`e}s
que $t(x)=s(y)$ pour tout $u\in U$ et tout $x,y\in \P X$.

Puis on adapte {\`a} la cat{\'e}gorie des flots la th{\'e}orie des paires NDR et
DR d'espaces topologiques \cite{MR35:970} de la fa{\c c}on suivante. La
notion de paire NDR d'espaces topologiques devient la notion de
\textit{S-cofibration} d{\'e}finie comme suit :

\begin{definition}\cite{flow} 
Un morphisme de flots $i:X\rightarrow Y$ est une S-cofibration si $i$
v{\'e}rifie la propri{\'e}t{\'e} de rel{\`e}vement des S-homotopies, i.e. pour tous
morphismes de flots $g:Y\rightarrow Z$ et $H:[0,1]\boxtimes
X\rightarrow Z$ tels que $H(0\boxtimes x)=g(i(x))$, il existe un
morphisme de flots $\overline{H}:[0,1]\boxtimes Y\rightarrow Z$ tel
que $\overline{H}(u\boxtimes i(x))=H(u\boxtimes x)$ et
$\overline{H}(0\boxtimes y)=g(y)$ pour tout $u\in [0,1]$, $x\in X$ et
$y\in Y$.
\end{definition}

La notion de paire DR d'espaces topologiques devient la notion de \textit{S-cofibration
acyclique} d{\'e}finie comme suit :

\begin{definition}\cite{flow} 
Une S-cofibration $i:X\rightarrow Y$ est acyclique si elle est une
{\'e}quivalence de S-homotopie. \end{definition}

Avec ces objets, on d{\'e}montre alors (entre autre) le th{\'e}or{\`e}me suivant :

\begin{theoreme}\cite{flow} 
Soit $(Y,B)$ une paire NDR d'espaces topologiques et soit
$i:A\rightarrow X$ une S-cofibration. Alors le morphisme de flots
canonique $(Y,B)\boxtimes i:Y\boxtimes A \sqcup_{B\boxtimes A}
B\boxtimes X\rightarrow Y\boxtimes X$ est une S-cofibration. Si de
plus $(Y,B)$ est une paire DR ou si $i$ est une S-cofibration
acyclique, alors la S-cofibration $(Y,B)\boxtimes i$ est
acyclique. \end{theoreme}

On a enfin besoin, comme dans le cas classique, de prouver que tout
morphisme de flots est S-homotope {\`a} une inclusion de flots. Cela se
fait en introduisant le "mapping cylindre" d'un morphisme de flots
comme suit :

\begin{definition}\cite{flow} 
Soit $f:X\rightarrow Y$ un morphisme de flots. Le mapping cylindre de
$f$ est le flot $I_f$ obtenu par le diagramme cocart{\'e}sien de $\dtop$
suivant :
\[
\xymatrix{\{1\}\boxtimes X \fr{}\fd{1\boxtimes x\mapsto f(x)}& [0,1]\boxtimes X\fd{}\\
Y \fr{} & \cocartesien I_f}
\]
\end{definition}

Muni de ces outils, on suit alors pas-{\`a}-pas la d{\'e}monstration du
th{\'e}or{\`e}me de Whitehead telle que expos{\'e}e dans le livre
\cite{MR80b:55001} en adaptant aux flots les th{\'e}or{\`e}mes de compression
ainsi que les diff{\'e}rents th{\'e}or{\`e}mes d'extension de morphismes ou
d'homotopies. On obtient m{\^e}me gratuitement le th{\'e}or{\`e}me suivant :

\begin{theoreme}\cite{flow} 
Soit $f:X\rightarrow Y$ un morphisme de flots o{\`u} $X$ et $Y$ sont deux
flots S-cofibrants. Alors $f$ est S-homotope {\`a} un morphisme de flots
$g:X\rightarrow Y$ tel que pour tout $n\geq 0$, $g(X^1_n)\subset
Y^1_n$. \end{theoreme}

Enfin on a :

\begin{theoreme}\cite{flow} Tout flot est faiblement S-homotope {\`a} un flot S-cofibrant.
Ce CW-complexe globulaire est n{\'e}cessairement unique {\`a} S-homotopie
pr{\`e}s d'apr{\`e}s le th{\'e}or{\`e}me pr{\'e}c{\'e}dent. \end{theoreme}

Si $X$ est un espace topologique, on sait qu'il existe une proc{\'e}dure
pour construire une famille de CW-complexes $(Y_n)_{n\geq 0}$ o{\`u} pour
tout $i\geq 0$ et tout point base, l'application canonique
$\pi_i(Y_n)\rightarrow \pi_i(X)$ est surjective et pour tout $0\leq
i\leq n$ et tout point base, l'application canonique
$\pi_i(Y_n)\rightarrow \pi_i(X)$ est injective. Dans cette proc{\'e}dure,
on passe de $Y_n$ {\`a} $Y_{n+1}$ en attachant des disques $D^{n+1}$ {\`a}
$Y_n$. La notation $\pi_i$ d{\'e}signe {\'e}videmment le $i$-i{\`e}me groupe
d'homotopie.

Une proc{\'e}dure exactement similaire fonctionne pour les flots. On part
d'un flot $X$ et on commence par associer {\`a} $\P X$ un espace $Y_0$ et
une application continue $Y_0\rightarrow \P X$ v{\'e}rifiant les m{\^e}mes
conditions que ci-dessus. On prend alors le flot libre $T_0$ engendr{\'e}
par $X^0$ et $Y_0$, ce qui ne change {\'e}videmment rien {\`a} la situation
pour les groupes d'homotopie. Puis on attache {\`a} $T_0$ des cellules
$Glob(D^1)$ le long de morphismes d'attachement $Glob(S^0)\rightarrow
T_0$ de fa{\c c}on {\`a} obtenir un morphisme de flots $T_1\rightarrow X$
induisant des bijections pour $\pi_0$ et $\pi_1$. L'unique diff{\'e}rence
avec ce qui se passe pour les espaces topologiques est qu'attacher une
cellule {\`a} un flot peut engendrer, du fait de la composition des
chemins, un attachement de plusieurs cellules {\`a} l'espace des chemins
sous-jacents.  Cela ne change rien pour les groupes d'homotopie. D'o{\`u}
le r{\'e}sultat en it{\'e}rant le proc{\'e}d{\'e}. Et comme corollaire, on obtient le

\begin{theoreme}\cite{flow} 
Le foncteur r{\'e}alisation cat{\'e}gorique induit une {\'e}quivalence de
cat{\'e}gories entre la cat{\'e}gorie des CW-complexes globulaires {\`a}
S-homotopie et {\`a} T-homotopie pr{\`e}s et la cat{\'e}gorie des flots {\`a}
S-homotopie faible et {\`a} T-homotopie pr{\`e}s. \end{theoreme}

Nous avons donc obtenu apr{\`e}s \cite{diCW} une deuxi{\`e}me fa{\c c}on de d{\'e}crire
les automates parall{\`e}les {\`a} dihomotopie pr{\`e}s. L'immense avantage de
cette seconde fa{\c c}on est que la cat{\'e}gorie des flots est compl{\`e}te et
cocompl{\`e}te, contrairement {\`a} la cat{\'e}gorie des CW-complexes
globulaires. De plus on a le

\begin{theoreme}\cite{flow}\label{mod} 
Il existe une structure mod{\`e}le sur la 
cat{\'e}gorie des flots dont les {\'e}quivalences faibles sont exactement 
les S-homotopies faibles. 
\end{theoreme}

Il existe une deuxi{\`e}me classe de morphismes de flots appel{\'e}s T-homotopie 
qu'on aimerait voir transformer en {\'e}quivalence faible. D'o{\`u} la

\begin{question} Existe-t-il une structure mod{\`e}le sur la cat{\'e}gorie des
flots dont les {\'e}quivalences faibles  contiendraient aussi les T-homotopies ? 
Cette structure mod{\`e}le pourrait {\^e}tre obtenue {\'e}ventuellement par localisation 
{\`a} gauche ou {\`a} droite au sens de \cite{ref_model2} de la structure mod{\`e}le du th{\'e}or{\`e}me~\ref{mod}. 
\end{question}

Il y a deux obstacles {\`a} surmonter pour r{\'e}pondre {\`a} cette question : 1) 
prouver l'existence d'une ou des deux localisations , 2) v{\'e}rifier que 
la ou les nouvelles structures mod{\`e}le n'identifient pas des automates 
parall{\`e}les ayant des propri{\'e}t{\'e}s diff{\'e}rentes.

Pour terminer, la cat{\'e}gorie des flots n'est pas cart{\'e}siennement ferm{\'e}e
mais elle poss{\`e}de quand m{\^e}me une structure mono{\"\i}dale ferm{\'e}e qui a une
signification informatique int{\'e}ressante \cite{flow}.

\providecommand{\bysame}{\leavevmode\hbox to3em{\hrulefill}\thinspace}
\providecommand{\MR}{\relax\ifhmode\unskip\space\fi MR }
\providecommand{\MRhref}[2]{%
  \href{http://www.ams.org/mathscinet-getitem?mr=#1}{#2}
}
\providecommand{\href}[2]{#2}


\begin{thebibliography}{1}

\bibitem{flow}
P.~Gaucher, \emph{{A Convenient Category for The Homotopy Theory of
  Concurrency}}, 2002, arXiv:math.AT/0201252.

\bibitem{note1}
P.~Gaucher, \emph{Automate parall{\`e}le {\`a} homotopie pr{\`e}s ({I})}, projet de note
  aux C.R.A.S., 2002.

\bibitem{diCW}
P.~Gaucher and E.~Goubault, \emph{{Topological Deformation of Higher
  Dimensional Automata}}, 2001, arXiv:math.AT/0107060, {\`a} para{\^i}tre dans
  Homology, Homotopy and Applications.

\bibitem{ref_model2}
P.~S. Hirschhorn, \emph{Localization of model categories}, available at
  http://www-math.mit.edu/\~{}psh/, October 2001.

\bibitem{Ref_wH}
L.~G. Lewis, \emph{The stable category and generalized thom spectra}, Ph.D.
  thesis, University of Chicago, 1978.

\bibitem{MR35:970}
N.~E. Steenrod, \emph{A convenient category of topological spaces}, Michigan
  Math. J. \textbf{14} (1967), 133--152.

\bibitem{MR80b:55001}
G.~W. Whitehead, \emph{Elements of homotopy theory}, Springer-Verlag, New York,
  1978.

\end{thebibliography}
\end{document}